\documentclass[12pt]{article}
\usepackage[utf8]{inputenc}
\usepackage{geometry,amsmath,amssymb,bm}
\usepackage[square, comma, sort&compress, numbers]{natbib}
\usepackage{natbib}
\usepackage{amsfonts}
\usepackage{latexsym}
\usepackage{graphicx}
\usepackage{amsmath}
\usepackage{color}
\usepackage{subfigure}
\usepackage{float}
\RequirePackage{jabbrv}
\newcommand{\beq}{\begin{equation}}
\newcommand{\bseqs}{\begin{subequations}}
\newcommand{\eseqs}{\end{subequations}}
\newcommand{\balign}{\begin{align}}
\newcommand{\ealign}{\end{align}}
\newcommand{\eeq}{\end{equation}}
\newcommand{\beql}{\begin{equation} \label}
\newcommand{\beqs}{\begin{eqnarray}}
\newcommand{\eeqs}{\end{eqnarray}}
\newcommand{\beas}{\begin{eqnarray*}}
\newcommand{\eeas}{\end{eqnarray*}}
\newcommand{\ber}{\begin{array}}
\newcommand{\eer}{\end{array}}
\newcommand{\becs}{\begin{cases}}
\newcommand{\eecs}{\end{cases}}
\newcommand{\leftm}{\left[\begin{array}}
\newcommand{\rightm}{\end{array}\right]}

\newcommand{\bfe}{{\mathbf e}}

\newcommand{\bfn}{{\mathbf n}}

\newcommand{\bfp}{{\mathbf p}}

\newcommand{\bfu}{{\mathbf u}}
\newcommand{\bfv}{{\mathbf v}}

\newcommand{\bfx}{{\mathbf x}}
\newcommand{\bfy}{{\mathbf y}}

\newcommand{\bfC}{{\mathbf C}}

\newcommand{\bfF}{{\mathbf F}}

\newcommand{\bfM}{{\mathbf M}}

\newcommand{\bfQ}{{\mathbf Q}}

\newcommand{\bfV}{{\mathbf V}}

\newcommand{\bbm}{\begin{bmatrix}}
\newcommand{\ebm}{\end{bmatrix}}

 \usepackage{mathptmx}
 \usepackage{color}
\newtheorem{theorem}{\bf Theorem}[section]

\newtheorem{proposition}{\bf Proposition}[section]
\begin{document}

\title{A constrained proof of the strong version of the Eshelby conjecture for the three-dimensional isotropic medium}
\author{T. Y. Yuan$^{1,2}$, K. F. Huang$^{3,*}$ and J. X. Wang$^{1,2,*}$
\bigskip
\\
\small{$^{1}$State Key Laboratory for Turbulence and Complex System, College of Engineering,}
\\
\small{Peking University, Beijing 100871, P.R. China}
\smallskip
\\
\small{$^{2}$CAPT-HEDPS, and IFSA Collaborative Innovation Center of MoE, College of Engineering,}
\\
\small{Peking University, Beijing 100871, P.R. China}
\smallskip
\\
\small{$^{3}$Department of Mechanics and Aerospace Engineering,}
\\
\small{Southern University of Science and Technology, Shenzhen, Guangdong 518055, P.R. China}
\\
}
%\listoftodos
\date{}
\maketitle

\begin{abstract}
Eshelby's seminal work on the ellipsoidal inclusion problem leads to the conjecture that
 the ellipsoid is the only inclusion possessing the uniformity property that a uniform eigenstrain is transformed into a uniform elastic strain.
 For the three-dimensional isotropic medium, the weak version of the Eshelby conjecture has been substantiated. The previous work of Ammari et al substantiates the strong version of the Eshelby conjecture for the cases when the three eigenvalues of the eigenstress are distinct or all the same, whereas the case where two of the eigenvalues of the eigenstress are identical and the other one is distinct remains a difficult problem. In this work, we study the latter case. To this end, firstly, we present and prove a necessary condition
 for a convex inclusion being capable of transforming a single uniform eigenstress  into a uniform elastic stress field.
 Since the necessary condition is not enough to determine the shape of the inclusion, secondly, we introduce a constraint that is concerned with the material parameters, and prove that
there exist combinations of the elastic tensors and uniform eigenstresses
such that only an ellipsoid can have the Eshelby uniformity property for these combinations simultaneously.
Finally, we provide a more specifically constrained proof of the conjecture by proving that
for the uniform strain fields constrained to that induced by an ellipsoid
from a set of specified uniform eigenstresses, the strong version of the Eshelby conjecture is true for a set of isotropic elastic tensors which are associated with the specified uniform eigenstresses. This work makes some progress towards the complete solution of the intriguing and longstanding Eshelby conjecture for three-dimensional isotropic media.
\end{abstract}

\thanks{
\begin{flushleft}
  \textbf{Subject Areas}: solid mechanics, applied mathematics, engineering\\
\end{flushleft}
\begin{flushleft}
\end{flushleft}
\textbf{Author for correspondence}: \\J. X. Wang: jxwang@pku.edu.cn;\quad K. F. Huang: huangkf@sustc.edu.cn
}

\section{Introduction}

Eshelby's  work~\cite{Eshelby1957,Eshelby1959} on the inclusion problem is essential to the development of the theories for the mechanical performance of heterogeneous materials. Mura introduced the concept of the eigenstrain, and conducted a series of studies via the Eshelby formalism~\cite{Mura1987}. Many researchers have investigated inclusion problems from diverse aspects. For instance, from the aspect of multi-physical fields, the generalized Eshelby formalisms for piezoelectric inclusions~\cite{Wang1992,Jiang1997} and inclusions governed by more general coupled-fields \cite{Li1998} have been developed; from the aspect of multi-scales,  the interface effect of nano-inclusions and nano-inhomogeneities has been studied~
\cite{Sharma2003,Sharma2004,Sharma2007,Nanoinc2005,Lim2006,Tian2007}, and the inclusion problem in the context of second gradient elasticity have  been investigated~\cite{Ma2018}. In addition, the classical Eshelby's inclusion theory for solids has been extended to investigations of liquid inclusions in soft materials \cite{Style2015,Ti2021}.

%, which can somehow be used to model the cancer cell influenced by the abnormal deposition and remodeling of extracellular matrix and thus helps to better understand the rate and direction of tumor cell migration .}

In the field of the inclusion problem, the most significant and fantastic phenomenon is the Eshelby uniformity property of the ellipsoid, which means that the uniform eigenstrain (eigenstress) prescribed in an ellipsoidal inclusion in an infinite medium induces  a uniform elastic strain (stress) field inside the inclusion.
In 1961, Eshelby \cite{Eshelby1961} conjectured that the ellipsoid uniquely possesses such marvelous uniformity property, which has long been a difficult problem to be proved or disproved.
According to previous researches \cite{Liu2008,Kang2008,Xu2009}, the Eshelby conjecture can be more specifically understood in two senses,  i.e. the weak and strong versions.
The weak version asserts that ``{\it an ellipsoid alone transforms all uniform eigenstrains (eigenstressess) into uniform elastic strain (stress) fields in it}'', and the strong version asserts that ``{\it no inclusion other than an ellipsoid transforms a single uniform eigenstrain (eigenstress) into a uniform elastic strain (stress) field in it}"~\cite{Yuan2021}. Note that there are also other ways of statements of the versions \cite{Liu2008,Kang2008,Xu2009}.
It is easily seen that the validity of the strong version leads to the validity of the weak version.

For the isotropic medium in three dimensions, in 2008, the proof of the weak version was fulfilled by Liu~\cite{Liu2008} with the utilization of the obstacle function on the basis of inclusion problem, and by Kang and Milton~\cite{Kang2008} with the utilization of the single-layer potential on the basis of the inhomogeneity problem. The conjecture is also valid for the three-dimensional conductivity problem \cite{Wang2014}.  For the isotropic medium in two dimensions, the proof of the Eshelby conjecture has been fulfilled~\cite{Sendeckyj1970,Ru1996,Vigdergauz2000}.
For anisotropic media, the weak and strong versions of the Eshelby conjecture have already been proved to be valid in two dimensions ~\cite{Xu2009}.
For the three-dimensional case, Yuan et al.~\cite{Yuan2021}  recently proved that the weak version for cubic, transversely isotropic, orthotropic, and monoclinic symmetries is valid, but there are counterexamples to the strong version.

However, the strong version of the conjecture for the three-dimensional isotropic case has not been fully resolved, though some progress has been made.
It was proved by Markenscoff  \cite{Markenscoff1997} in 1997 that the inclusion that possesses the Eshelby uniformity property can not have any planus surface. In 1998, Lubarda and Markenscoff~\cite{Lubarda1998} further drew the conclusion that the surface of the inclusion that possesses the Eshelby uniformity property needs to satisfy some particular conditions.  Then Markenscoff~\cite{Markenscoff1998} also proved that the only way to assure the Eshelby uniformity property is the infinitesimal perturbation of the ellipsoid into another ellipsoid.

%And Liu~\cite{Liu2008}, and Kang and Milton~\cite{Kang2008} both substantiated the weak version by investigating into a related Newtonian %potential problem that corresponds to the inclusion domain.
In 2010, a further step towards the proof of the strong version was made by Ammari et al.~\cite{Ammari2010}.
By categorizing the induced strain field based on its eigenvalues, Ammari et al.~\cite{Ammari2010} proved that when the three eigenvalues of the elastic strain field induced by the remote loading are either all the same or all distinct, the ellipsoidal inhomogeneity uniquely possesses the Eshelby uniformity property, which is stronger than the weak version and close to the strong version.

It is noted that the result given by Ammari et al.~\cite{Ammari2010} for the inhomogeneity problem in the isotropic medium can be extended to prove that the ellipsoidal inclusion uniquely possesses the Eshelby uniformity property, when the eigenvalues of the eigenstress are either all distinct or all identical, that is, the strong version of the Eshelby conjecture is true for these two special cases. In this work, we consider the remaining case where two of the eigenvalues of the eigenstress are identical, and the other one is distinct. In this case, we show and prove three theorems.
The first theorem constitutes  a necessary condition for a convex inclusion to possess the Eshelby uniformity property. However, the necessary condition alone can not help us to determine the shape of the inclusion. Therefore, in the second theorem,
we bring in an additional constraint associated with the material parameters, and find that the inclusion can only be an ellipsoid owing to the necessary condition under the constraint that the inclusion possesses the Eshelby uniformity property for two isotropic media whose elastic tensors are linearly independent.  Thus, the second theorem provides a constrained proof of the strong version.
Further, we evaluate the elastic field induced by an arbitrary inclusion with the Eshelby uniformity property, and find that the uniform strain field induced by an inclusion of any non-ellipsoidal shape, if there is one, can not be equal to that induced by an ellipsoid for the same uniform eigenstress, which constitutes the third theorem.  Thus the third theorem provides an alternative constrained proof of the strong version from another viewpoint when the induced strain field is constrained to be equal to that induced by an ellipsoid for the same uniform eigenstress.
All together, this work makes some progress towards the complete resolution of the strong version of the Eshelby conjecture for the three-dimensional
isotropic medium.

\section{Basic equations}

Let $\Omega\subset \mathbb{R}^{3}$ denote the inclusion domain, which is a one-component connected bounded domain with a Lipschitz boundary.
The equilibrium equation  for Eshelby's inclusion problem in the infinite elastic homogenous isotropic medium is
\begin{align}\label{equilibrium}
 \boldsymbol{\nabla}\cdot[\bfC:(\nabla\otimes\bfu)-\chi_\Omega\boldsymbol{\sigma}^*]={\bf 0} \;\;\mathrm{in} \;\;\;\mathbb{R}^{3},
\end{align}
where $\bfu$ denotes the displacement field, which is a vector; $\nabla$ denotes the gradient operator; $\boldsymbol{\sigma}^*$ denotes a uniform eigenstress, which is a symmetric second-order tensor; $\chi_\Omega$ denotes the indicator function with respect to $\Omega$; and
\begin{align}\label{isoelastictensor}
\bfC=[\lambda\delta_{ij}\delta_{kl}+\mu(\delta_{ik}\delta_{jl}+\delta_{il}\delta_{jk})]\;\bfe_i\otimes\bfe_j\otimes\bfe_k\otimes\bfe_l
\end{align}
denotes the fourth-order elastic tensor for the isotropic meidum with $\delta_{ij}$ being the Kronecker delta, $\lambda$ and $\mu$ being Lam\'{e} parameters that satisfy
 \begin{align}\label{Lam}
\begin{split}
\mu>0, \;\;3\lambda+2\mu>0,
\end{split}
\end{align}
and $\{\bfe_1,\bfe_2,\bfe_3\}$ being the bases of a Cartesian coordinate system $\bfx=(x_1,x_2,x_3)$. Unless otherwise stated, the summation convention on repeated indices will always be stipulated.

By the Fourier analysis \cite{Liu2008,Yuan2021}, the solution to \eqref{equilibrium} is
\begin{align}\label{displacement}
u_p(\bfx)=-\frac{\mathrm{i}}{(2\pi)^3}\int_{\mathbb{R}^3} L_{pq}(\boldsymbol{\xi}){\sigma}^*_{qj}\xi_j\int_{\Omega}e^{\mathrm{i}(\bfx-\bfy)\cdot\boldsymbol{\xi}}\mathrm{d}\bfy\mathrm{d}\boldsymbol{\xi},
\end{align}
where $\mathrm{i}=\sqrt{-1}$ denotes the imaginary unit, and
\begin{align}\label{L}
L_{pq}C_{qlmn}\xi_l\xi_n=\delta_{pm}.
\end{align}
By substituting \eqref{isoelastictensor} into \eqref{L}, we gain
\begin{align}\label{Green iso}
L_{pq}(\boldsymbol{\xi})=\frac{1}{\mu |\boldsymbol{\xi}|^2}\delta_{pq}-\frac{\mu+\lambda}{\mu(2\mu+\lambda)}\frac{\xi_p\xi_q}{ |\boldsymbol{\xi}|^4}.
\end{align}
Owing to the isotropy of the elastic tensor $\bfC$, which means $\bfC$ possesses the same expression \eqref{isoelastictensor} in all Cartesian coordinate systems, we see that \eqref{Green iso} and thus \eqref{displacement} do not vary with the rotation of the Cartesian coordinate system. Thus to simplify the derivations in the sequel, we let the axes of the Cartesian coordinate system  $\bfx=(x_1,x_2,x_3)$  coincide with the three mutually orthogonal unit eigenvectors of the eigenstress $\boldsymbol{\sigma}^*$. Then in such a coordinate system, the off-diagonal elements of $\boldsymbol{\sigma}^*$ are zero, and the diagonal elements $\sigma^*_{11}, \sigma^*_{22}, \sigma^*_{33}$ are  the three eigenvalues of $\boldsymbol{\sigma}^*$.
%can be expressed as
%\begin{align}\label{eigenstressineigenspace}
%\boldsymbol{\sigma}^*=\sigma^*_{JJ}\bfe_j\otimes\bfe_j,
%\end{align}
%where $\sigma^*_{JJ}\;(J=1,2,3)$ denote the three eigenvalues of the eigenstress $\boldsymbol{\sigma}^*$, and the capital letter $J$
%does not obey the summation convention but takes the same index as that of the corresponding  letter $j$ in \eqref{eigenstressineigenspace}.

%The formula \eqref{displacement} along with \eqref{eigenstressineigenspace} allows us to express  $\bfu$ only by the eigenvalues $\sigma^*_{JJ}\;(J=1,2,3)$ of the eigenstress $\boldsymbol{\sigma}^*$ and the shape of the inclusion $\Omega$.
In this coordinate system, the inclusion $\Omega$ that possesses the Eshelby uniformity property must yield
\begin{align}\label{finaluineigenspace}
\begin{split}
\frac{\partial u_p(\bfx)}{\partial x_l}\;&=\frac{1}{(2\pi)^3}\int_{\mathbb{R}^3} L_{pq}(\boldsymbol{\xi}){\sigma}^*_{qj}\xi_j\xi_l\int_{\Omega}e^{\mathrm{i}(\bfx-\bfy)\cdot\boldsymbol{\xi}}\mathrm{d}\bfy\mathrm{d}\boldsymbol{\xi}=\;\;\mathrm{constant},
\end{split}
\end{align}
%where ${\sigma}^*_{qj}=0$ for $q\neq j$.
%where ${\sigma}^*_{qj}$ denotes a single uniform eigenstress shown in \eqref{eigenstressineigenspace}.
Note that the left-hand side of \eqref{finaluineigenspace} corresponds to the total strain induced by $\Omega$. Since the total strain is the sum of the elastic strain and the uniform eigenstrain, thus the uniformity of the total strain is equivalent to the uniformity of the elastic strain.
Thus based on \eqref{finaluineigenspace}, the completion of the proof of the strong version of the Eshelby conjecture for the three-dimensional isotropic medium can be achieved by classifying the three eigenvalues $\sigma^*_{11}, \sigma^*_{22}, \sigma^*_{33}$ of the eigenstress into three cases: they are all identical; they are all distinct; and two of them are identical and the other one is distinct, and then proving $\Omega$ that leads to \eqref{finaluineigenspace} must be ellipsoidal for these cases.

As is mentioned before, the strong version for the eigenstress possessing either all identical or all distinct eigenvalues can be proved by extending the result of Ammari et al. \cite{Ammari2010} to the inclusion problem via the same mathematical manipulations in Section 3 of their paper.  And we also provide an alternative proof for this case by using the Fourier analysis in the Appendix.
Therefore, only the case where two eigenvalues of the eigenstress are identical and the other one is distinct will be studied in detail in the following sections.
 Note that the classification of the eigenstress is equivalent to the classification of the eigenstrain due to the isotropy of the elastic tensor.

\section{A constrained proof of the strong version for the eigenstress possessing only two identical eigenvalues}
We define
\begin{align}\label{gsimplification}
g(\bfx,\boldsymbol{\xi}):=\int_{\Omega}e^{\mathrm{i}(\bfx-\bfy)\cdot\boldsymbol{\xi}}\mathrm{d}\bfy
\end{align}
and let $\sigma^*_{11}=\sigma^*_{22}=k_1,\;\sigma^*_{33}=k_3$ with $k_1\neq  k_3$.

Since the inclusion $\Omega$ with the Eshelby uniformity property satisfies \eqref{finaluineigenspace}, then
by substituting $\sigma^*_{11}=\sigma^*_{22}=k_1,\;\sigma^*_{33}=k_3$ along with \eqref{Green iso} and  \eqref{gsimplification} into  \eqref{finaluineigenspace}, we obtain

\begin{eqnarray}\label{twoidentical}
\begin{split}
\forall j=1,2,3\quad&\int_{\mathbb{R}^3} \xi_1\xi_j \left(\frac{\mu k_1}{|\boldsymbol{\xi}|^2}+\frac{(\lambda+\mu)(k_1-k_3)\xi_3^2}{|\boldsymbol{\xi}|^4}\right)g(\bfx,\boldsymbol{\xi})\mathrm{d}\boldsymbol{\xi}=\mathrm{constant},\\
&\int_{\mathbb{R}^3} \xi_2\xi_j \left(\frac{\mu k_1}{|\boldsymbol{\xi}|^2}+\frac{(\lambda+\mu)(k_1-k_3)\xi_3^2}{|\boldsymbol{\xi}|^4}\right)g(\bfx,\boldsymbol{\xi})\mathrm{d}\boldsymbol{\xi}=\mathrm{constant},\\
  &\int_{\mathbb{R}^3} \xi_3\xi_j\left(\frac{ k_3(\lambda+2\mu)-k_1(\lambda+\mu)}{|\boldsymbol{\xi}|^2}+\frac{(\lambda+\mu)(k_1-k_3)\xi_3^2}{|\boldsymbol{\xi}|^4}\right)g(\bfx,\boldsymbol{\xi})\mathrm{d}\boldsymbol{\xi}=\mathrm{constant}\quad\bfx\in \Omega.\\
 \end{split}
\end{eqnarray}

 It follows from \eqref{twoidentical} that

\begin{eqnarray}\label{disfourier}
\begin{split}
\quad&\frac{\partial}{\partial x_1}\int_{\mathbb{R}^3} \left(\frac{\mu k_1}{|\boldsymbol{\xi}|^2}+\frac{(\lambda+\mu)(k_1-k_3)\xi_3^2}{|\boldsymbol{\xi}|^4}\right)g(\bfx,\boldsymbol{\xi})\mathrm{d}\boldsymbol{\xi}=\mathrm{Linear},\\
\quad &\frac{\partial}{\partial x_2}\int_{\mathbb{R}^3} \left(\frac{\mu k_1}{|\boldsymbol{\xi}|^2}+\frac{(\lambda+\mu)(k_1-k_3)\xi_3^2}{|\boldsymbol{\xi}|^4}\right)g(\bfx,\boldsymbol{\xi})\mathrm{d}\boldsymbol{\xi}=\mathrm{Linear},\\
 \quad &\frac{\partial}{\partial x_3}\int_{\mathbb{R}^3}\left(\frac{ k_3(\lambda+2\mu)-k_1(\lambda+\mu)}{|\boldsymbol{\xi}|^2}+\frac{(\lambda+\mu)(k_1-k_3)\xi_3^2}{|\boldsymbol{\xi}|^4}\right)g(\bfx,\boldsymbol{\xi})\mathrm{d}\boldsymbol{\xi}=\mathrm{Linear}\quad\bfx\in \Omega,\\
 \end{split}
 \end{eqnarray}

\begin{flushleft}
 which finally leads to
\end{flushleft}

\begin{align}\label{Ntilde2}
\begin{split}
\quad &\int_{\mathbb{R}^3} \left(\frac{\mu k_1}{|\boldsymbol{\xi}|^2}+\frac{(\lambda+\mu)(k_1-k_3)\xi_3^2}{|\boldsymbol{\xi}|^4}\right)g(\bfx,\boldsymbol{\xi})\mathrm{d}\boldsymbol{\xi}=q_1(\bfx)+\varphi_1(x_3),\\
 \quad &\int_{\mathbb{R}^3}\left(\frac{ k_3(\lambda+2\mu)-k_1(\lambda+\mu)}{|\boldsymbol{\xi}|^2}+\frac{(\lambda+\mu)(k_1-k_3)\xi_3^2}{|\boldsymbol{\xi}|^4}\right)g(\bfx,\boldsymbol{\xi})\mathrm{d}\boldsymbol{\xi}=q_2(\bfx)+\varphi_2(x_1,x_2)\quad\bfx\in \Omega,\\
 \end{split}
\end{align}

\begin{flushleft}
  where $q_1(\bfx)$ and $q_2(\bfx)$ denote two quadratic functions, and $\varphi_1(x_3)$ and $\varphi_2(x_1,x_2)$ denote two unknown functions.
\end{flushleft}

To continue the analysis, we introduce two potentials.
The first one is the well-known Newtonian potential $N_\Omega(\bfx)$ of the inclusion domain $\Omega$, which can be expressed as
\begin{align}\label{NewtonianFourier}
N_\Omega(\bfx)=-\frac{1}{(2\pi)^3}\int_{\mathbb{R}^3} \frac{g(\bfx,\boldsymbol{\xi})}{|\boldsymbol{\xi}|^2}\mathrm{d}\boldsymbol{\xi}.
\end{align}
Note that the Newtonian potential $N_\Omega(\bfx)$ in \eqref{NewtonianFourier} is exactly the first term of the integral on the left-hand side of \eqref{Ntilde2}. Besides, we know that the Newtonian potential $N_\Omega(\bfx)$ satisfies
\begin{align}\label{definition N}
\begin{split}
\left\{ \begin{array}{*{20}{l}}
&{\Delta N_\Omega(\bfx)=\chi_\Omega \quad\quad\;\;\;\bfx\in \mathbb{R}^3,}\\
&{N_\Omega(\bfx)= \mathcal{O}\left(\frac{1}{|\bfx|}\right)\quad \text{as} \;|\bfx|\to+\infty,}
\end{array}  \right.
\end{split}
\end{align}
where $\Delta$ denotes the Laplacian operator, and $\mathcal{O}(\cdot)$ denotes the order of magnitude.  \eqref{definition N} admits a unique solution, i.e.,
\begin{align}\label{eq:TheoremNP}
 N_\Omega(\bfx)=-\frac{1}{4\pi}\int_{\Omega} \frac{1}{{|\bfx-\bfy|}} \mathrm{d}\bfy,
 \end{align}
which  provides an explicit expression of $N_\Omega(\bfx)$.

The second one is a bi-harmonic potential $H$, which is defined as
\begin{align}\label{Hex}
H(\bfx):= \frac{1}{(2\pi)^3}\int_{\mathbb{R}^3} \frac{g(\bfx,\boldsymbol{\xi})}{|\boldsymbol{\xi}|^4}\mathrm{d}\boldsymbol{\xi}.
\end{align}
By substituting \eqref{Hex} into \eqref{NewtonianFourier}, it is easy to verify that
\begin{align}\label{definition HN}
\Delta H(\bfx)= N_\Omega(\bfx);
\end{align}
thus by $\eqref{definition N}_1$,
\begin{align}\label{definition H}
\Delta^2 H(\bfx)=\chi_\Omega\quad\bfx\in\mathbb{R}^3,
\end{align}
which indicates $H(\bfx)$ is bi-harmonic in $\mathbb{R}^3\setminus\Omega$.

%which means that $H(\bfx)$ at infinity is proportional to the length $|\bfx|$ of $\bfx$ and not associated with the direction of $\bfx$.}

%Note that the existence of the Fourier transformation of $H(\bfx)$ and thus the Fourier form \eqref{Hex} of $H(\bfx)$ call for some properties of $H(\bfx)$, such as the integrability and the boundary condition, which are not listed here.

It could be derived from \eqref{Hex} that
\begin{align}\label{Htwoderivative}
\frac{\partial^2 H(\bfx)}{\partial x_3^2}= -\frac{1}{(2\pi)^3}\int_{\mathbb{R}^3} \frac{\xi_3^2}{|\boldsymbol{\xi}|^4}g(\bfx,\boldsymbol{\xi})\mathrm{d}\boldsymbol{\xi},
\end{align}
which is exactly the second term of the integral on the left-hand side of \eqref{Ntilde2}. Therefore, substitution of \eqref{NewtonianFourier} and \eqref{Htwoderivative} into \eqref{Ntilde2} yields
\begin{align}\label{NandH}
\begin{split}
 &\alpha N_\Omega(\bfx)+\frac{\partial^2 H(\bfx)}{\partial x_3^2}=q_1(\bfx)+\varphi_1(x_3),\\
 &\beta N_\Omega(\bfx)+\frac{\partial^2 H(\bfx)}{\partial x_3^2}=q_2(\bfx)+\varphi_2(x_1,x_2)\quad\bfx\in \Omega,
 \end{split}
\end{align}
where $\alpha$ and $\beta$ are two real constants defined by
\begin{align}\label{alphabeta}
\begin{split}
&\alpha:=-\frac{\mu k_1 }{(\lambda+\mu)(k_1-k_3)};\\
& \beta:=-\frac{[k_3(\lambda+2\mu)-k_1(\lambda+\mu)]}{(\lambda+\mu)(k_1-k_3)}.
\end{split}
\end{align}
 Note that, here,  the values of $q_i$ and $\varphi_i$ in \eqref{NandH} are equal to  $\frac{1}{(2\pi)^3}$ times  $q_i$ and $\varphi_i$ in \eqref{Ntilde2}.
Since $k_1\neq k_3$, it is easy to verify that for any combination of the Lam\'{e} parameters $\lambda,\mu$ that satisfy \eqref{Lam}, \eqref{alphabeta} is always valid, and $\alpha\neq \beta$.

To determine the shape of $\Omega$, we proceed to analyze \eqref{NandH}.
Previously, we have shown the expression \eqref{eq:TheoremNP} of the Newtonian potential $N_\Omega(\bfx)$.  By substituting \eqref{eq:TheoremNP} into \eqref{definition HN}, we can derive
\begin{align}\label{Hexpression0}
H(\bfx)=-\frac{1}{8\pi}\int_{\Omega}|\bfx-\bfy|\mathrm{d}\bfy+H^*(\bfx),
\end{align}
where $H^*(\bfx)$  is an unknown harmonic function that satisfies
\begin{align}\label{H*}
\begin{split}
\Delta H^*(\bfx)=0 \quad\bfx\in \mathbb{R}^3.
\end{split}
\end{align}
Then we are going to determine $H^*(\bfx)$.

By combining \eqref{definition HN} with $\eqref{definition N}_2$, we see
\begin{align}\label{Hinftyproperty1}
H(\bfx)= \mathcal{O}\left(|\bfx|\right)\quad \text{as} \;\;|\bfx|\to+\infty.
\end{align}
And by \eqref{Hex} and \eqref{Hinftyproperty1}, it is easy to verify that
\begin{align}\label{Hinftyproperty2}
\begin{split}
H(\bfx)\;\propto \;|\bfx| \quad \text{as} \;\;|\bfx|\to+\infty.\\
\end{split}
\end{align}
Further, by substituting \eqref{Hexpression0} along with \eqref{H*} into \eqref{Hinftyproperty2}, we see $H^*(\bfx)$ at infinity must satisfy either
\begin{align}\label{H*bounday}
\begin{split}
H^*(\bfx)\propto |\bfx| \quad as\;\;|\bfx|\to+\infty,
\end{split}
\end{align} or
\begin{align}\label{H*bounday1}
\begin{split}
H^*(\bfx)\to 0 \quad as\;\;|\bfx|\to+\infty.
\end{split}
\end{align}
%which indicates that $H^*(\bfx)$ at infinity is unrelated to the direction of $\bfx$.
%Since a harmonic function is analytic~\cite{Han2012}, in the neighbourhood $U(\bfx_0)$ of any point $\bfx_0$, $H^*(\bfx)$ can be expressed as
%\begin{align}
%\begin{split}
%H^*(\bfx)=C_0+\bfA_0\cdot\bfx+W_0(\bfx)\quad\bfx\in U(\bfx_0),
%\end{split}
%\end{align}
% where $C_0$ denotes a constant;  $\bfA_0$ denotes a constant second-order matrix; and $\hat{\mathcal{O}}( \cdot)$ denotes the sum of the terms whose order is larger than 2.  $C_0$, $\bfA_0$ and the coefficients of the terms in $W_0(|\bfx|)$ are all determined by $\bfx_0$.

Take a ball $B_r=\{\bfx\;|\;|\bfx|\leq r,\;r>0\}$. Due to the mean-value property of the harmonic function, we have
\begin{align}\label{Meanproperty}
H^*({\bf 0})=\frac{1}{4\pi r^2}\int_{\partial B_r} H^*(\bfx) \mathrm{d}\bfx.
\end{align}
We then  consider the case when $r$ tends to infinity in \eqref{Meanproperty}.  If \eqref{H*bounday} holds, by substituting \eqref{H*bounday} into \eqref{Meanproperty} with $r\to +\infty$, we obtain
\begin{align}
H^*({\bf 0})\propto r\quad \text{where}\;\;r\to+\infty.
\end{align}
The above result is impossible since $H({\bf 0})$ is a finite constant. Hence \eqref{H*bounday} can not be valid, and thus $H^*(\bfx)$ at infinity should satisfy \eqref{H*bounday1}.

%It is easy to verify that \eqref{Hexpression0} is the expression of  $H(\bfx)$ in the real space. To this end, we assume that there is an additional term $H^*(\bfx)$ of the bi-harmonic potential $H(\bfx)$. Since the Newtonian potential $N_\Omega(\bfx)$ on the right-hand side of \eqref{definition HN} is uniquely determined by $\Omega$, thus it can be derived from \eqref{definition HN} that

%Since the Fourier form \eqref{Hex} of $H(\bfx)$ is uniquely determined by $\Omega$; thus for a given $\Omega$, $H^*(\bfx)$ should not contribute to the boundary value of $H(\bfx)$. Hence
%\begin{align}\label{H**}
%\begin{split}
%H^*(\bfx)= 0\quad\quad|\bfx|\to\infty.
%\end{split}
%\end{align}
Equations~\eqref{H*} and \eqref{H*bounday1} signify that $H^*(\bfx)$ is a harmonic function in $\mathbb{R}^3$, and $H^*(\bfx)$  will tend to zero at infinity; thus due to the analyticity of $H^*(\bfx)$ in $\mathbb{R}^3$, we see $H^*(\bfx)$ is bounded in $\mathbb{R}^3$. Then owing to the boundary condition \eqref{H*bounday1} and the Liouville theorem,  which stipulates that ``{\it any harmonic function in $\mathbb{R}^3$ bounded from above or
below is constant}"~\cite{Han2012},  we conclude that
\begin{align}\label{Hfinalexpression}
H^*(\bfx)\equiv 0.
\end{align}

It can be derived from \eqref{Hexpression0} along with \eqref{Hfinalexpression} that
\begin{align}\label{Hexpression}
\frac{\partial^2H(\bfx)}{\partial x_3^2}=-\frac{1}{8\pi}\int_{\Omega}\left[\frac{1}{|\bfx-\bfy|}-\frac{(x_3-y_3)^2}{|\bfx-\bfy|^3}\right]\mathrm{d}\bfy.
\end{align}
Given this, we define
\begin{align}\label{Ntildeexpression}
\tilde{N}_{\Omega}(\bfx):=-\frac{1}{4\pi}\int_{\Omega }\frac{(x_3-y_3)^2}{|\bfx-\bfy|^3}\mathrm{d}\bfy.
\end{align}
Then by substituting \eqref{eq:TheoremNP}, \eqref{Hexpression} and \eqref{Ntildeexpression} into \eqref{NandH}, we obtain
\begin{align}\label{NandNtilde}
\begin{split}
 \gamma N_\Omega(\bfx)+\tilde{N}_{\Omega}(\bfx)&=q_1(\bfx)+\varphi_1(x_3),\\
\eta N_\Omega(\bfx)+\tilde{N}_{\Omega}(\bfx)&=q_2(\bfx)+\varphi_2(x_1,x_2)\quad\bfx\in \Omega,
 \end{split}
\end{align}
where $N_\Omega(\bfx)$ and $\tilde{N}_{\Omega}(\bfx)$ both have explicit expressions, i.e., \eqref{eq:TheoremNP} and \eqref{Ntildeexpression}, respectively, and
  \begin{align}\label{11.18}
\begin{split}
&{{\gamma }}:=\frac{[k_3(\lambda+\mu)-k_1(\lambda+3\mu)]}{(\lambda+\mu)(k_1-k_3)},\\
& {{\eta }} :=\frac{[k_3(\lambda+3\mu)-k_1(\lambda+\mu)]}{(\lambda+\mu)(k_1-k_3)}
\end{split}
\end{align}
are real constants that are determined by the elastic tensors $\bfC$ (Lam\'{e} parameters $\lambda,\mu$) and the eigenvalues $k_1,k_3$ of the eigenstress $\boldsymbol{\sigma}^*$.

Recently, Yuan et al.~\cite{Yuan2021} demonstrate that the material parameters can serve as an important factor in the study of the Eshelby conjecture. Likewise, we also consider the influence of the material parameters in this work. To prove the strong version of the Eshelby conjecture for the eigenstress possessing only two identical eigenvalues is to prove the following proposition:

\begin{proposition}\label{pro1}
Let $\Omega\subset \mathbb{R}^{3}$ be a one-component connected bounded open  domain with a Lipschitz boundary. Equation \eqref{NandNtilde} holds
for $(\boldsymbol{\sigma}^*,\bfC)$, where $\boldsymbol{\sigma}^*$ is a given  uniform eigenstress that possesses two identical eigenvalues $k_1$ and a distinct eigenvalue $k_3$, and $\bfC$ is an isotropic elastic tensor,  if and only if $\Omega$ is of ellipsoidal shape.
\end{proposition}
%Note that Proposition~\ref{pro1} sustains the strong version of the Eshelby conjecture for the eigenstress only possessing two identical %eigenvalues. Since we have already gained the proof of the strong version for the other two cases, thus the proof of Proposition~\ref{pro1} %consequently completes the proof of the strong version of the Eshelby conjecture.

To prove or disprove Proposition~\ref{pro1}, we investigate \eqref{NandNtilde}.
Firstly, we derive a necessary condition for a convex inclusion possessing the Eshelby uniformity property from \eqref{NandNtilde}, which is stated in the following theorem:

\begin{theorem}\label{T0.2}
Let $\Omega\subset \mathbb{R}^{3}$ be a one-component connected bounded open convex domain with a Lipschitz boundary. If equation \eqref{NandNtilde} holds
for $(\boldsymbol{\sigma}^*,\bfC)$ where $\boldsymbol{\sigma}^*$ is a given uniform eigenstress that possesses two identical eigenvalues $k_1$ and a distinct eigenvalue $k_3$, and $\bfC$ is an isotropic elastic tensor, then in the Cartesian coordinate system $\bfx=(x_1,x_2,x_3)$ defined by the eigenvectors of the eigenstress $\boldsymbol{\sigma}^*$, there must exist an ellipsoid $E \subset \Omega$ that satisfies
\begin{align}\label{Ncorrelation}
N_{\Omega}(\bfx)-N_{E}(\bfx)= {\varphi}(x_1,x_2), \;\; \bfx\in E,
\end{align}
where $\varphi(x_1,x_2)$ is an unknown function that only relies on two spatial coordinates, and $N_{E}(\bfx)$  is the Newtonian potential of the ellipsoid $E$, whose expression is \eqref{eq:TheoremNP} with $\Omega$ replaced by $E$.
\end{theorem}

Theorem \ref{T0.2} implies that the Newtonian potential of a convex $\Omega$ that possesses the Eshelby uniformity property for the eigenstress possessing only two identical eigenvalues is necessarily correlated with the Newtonian potential of an ellipsoid $E$ via \eqref{Ncorrelation}.
If $\Omega$ is ellipsoidal, \eqref{Ncorrelation} is automatically satisfied.
%owing to the quadratic form of the Newtonian potential of ellipsoids \cite{Ferrers,Dyson}.
Thus, the significance of Theorem \ref{T0.2} is that a non-ellipsoidal convex inclusion that does not satisfy \eqref{Ncorrelation} is excluded from the set of inclusions that possess the Eshelby uniformity property.

Secondly, we impose an extra constraint concerned with the material parameters to make \eqref{NandNtilde} the sufficient and necessary condition  for $\Omega$ to be the ellipsoidal shape, which is stated in the following theorem:

\begin{theorem}\label{new}
Let $\Omega\subset \mathbb{R}^{3}$ be a one-component connected bounded open domain with a Lipschitz boundary. There exist combinations $( \boldsymbol{\sigma^*},\bfC^{(1)})$ and $( \boldsymbol{\sigma^*},\bfC^{(2)})$,  where $\boldsymbol{\sigma}^*$ is a given  uniform eigenstress that possesses two identical eigenvalues $k_1$ and a distinct eigenvalue $k_3$, and $\bfC^{(1)}$ and $\bfC^{(2)}$
are {two  linearly independent isotropic elastic tensors}, such that  equation \eqref{NandNtilde} holds for
$( \boldsymbol{\sigma^*},\bfC^{(1)})$ and $( \boldsymbol{\sigma^*},\bfC^{(2)})$ simultaneously,
if and only if $\Omega$ is of ellipsoidal shape.
\end{theorem}

By comparing Theorem \ref{new} with Proposition~\ref{pro1}, it is straightforward to see that the condition that \eqref{NandNtilde} holds for $( \boldsymbol{\sigma^*},\bfC^{(1)})$ and $( \boldsymbol{\sigma^*},\bfC^{(2)})$ simultaneously is stronger than the condition that \eqref{NandNtilde} holds for  $( \boldsymbol{\sigma^*},\bfC)$, which reveals extra constraints here. Thus Theorem \ref{new} provides a constrained proof of Proposition~\ref{pro1} and thus the strong version of the Eshelby conjecture for the eigentress possessing only two identical eigenvalues.

It is seen that the above theorems do not exclude the existence of non-ellipsoidal inclusions that can satisfy the Eshelby uniformity property (otherwise the strong version is ultimately proved ). Then, if there exist non-ellipsoidal inclusions that satisfy the Eshelby uniformity property,
what are the (uniform) strain fields in them like? We explore the answer to this interesting question by the following theorem:

\begin{theorem}\label{T0.3}
Let $\Omega\subset \mathbb{R}^{3}$ be a one-component connected bounded open domain with a Lipschitz boundary. There exists an isotropic elastic tensor $\bfC$ and an ellipsoid $E$ with $E\supset\Omega$, such that the elastic strain field induced by $\Omega$ is equal to that induced by $E$ within $\Omega$ for a given uniform eigenstress $\boldsymbol{\sigma}^*$  that possesses two identical eigenvalues $k_1$ and a distinct eigenvalue $k_3$ with the same sign as that of $k_1$,  if and only if $\Omega$ is of ellipsoidal shape.
\end{theorem}

Theorem \ref{T0.3} implies that for the uniform strain fields constrained to be equal to that induced by an ellipsoid from the specified $k_1$ and $k_3$, the strong version of the
Eshelby conjecture is true for a set of isotropic elastic tensors $\bfC$, which, as we will show below,
are determined by the eigenvalues $k_1$ and $k_3$.

 Theorem \ref{T0.3} can be expressed mathematically as follows.
We define a mapping $F$ that maps the Cartesian product of the set $\{\boldsymbol{\sigma ^*}\}$ of the uniform eigenstress, the set $\{\bfC\}$  of the isotropic elastic tensor and the set $\{\Omega\}$  of the configuration of the inclusion into the set $\{\boldsymbol{\varepsilon}^e\}$ of the induced elastic strain  inside $\Omega$ for the eigenstrain problem, i.e.,
\beqs\label{mapping}
\begin{split}
F:&\{\boldsymbol{\sigma ^*}\}\times \{\bfC\}\times \{\Omega\} \rightarrow\{\boldsymbol{\varepsilon}^e\},\;(\boldsymbol{\sigma}^*,\;\bfC,\;\Omega) \mapsto F(\boldsymbol{\sigma}^*,\;\bfC,\;\Omega) \subset {\mathbb{R}^{3\times3}}.
\end{split}
\eeqs

 %According to the above derivations, we have
%\beqs\label{mapping1}
%\begin{split}
%\varepsilon_{pl} &\equiv F(\boldsymbol{\sigma ^*},\;\bfC,\;\Omega)\\
%=&\frac{\int_{\mathbb{R}^3} \Big(L_{pq}(\boldsymbol{\xi}){\sigma}^*_{qj}\xi_j\xi_l+L_{lq}(\boldsymbol{\xi}){\sigma}^*_{qj}\xi_j\xi_p\Big)\frac{g(\bfx,\boldsymbol{\xi})}{2}\mathrm{d}\boldsymbol{\xi}}{(2\pi)^3}\\
%&-C_{plmn}\sigma ^*_{mn},
%\end{split}
%\eeqs
%which provides an explicit expression of $F$.
As is known from the previous researches \cite{Eshelby1957,Eshelby1959,Mura1987},
\beqs
\begin{split}
&\forall\;(\boldsymbol{\sigma^*},\;\bfC,\;E)\in \{\boldsymbol{\sigma ^*}\}\times\{\bfC\}\times\{E\},\;F(\boldsymbol{\sigma^*},\;\bfC,\;E)=\mathrm{constant},
 \end{split}
\eeqs
where $E$ denotes an ellipsoid, and $\{E\}$ denotes the set of ellipsoids, which is a subset of $\{\Omega\}$.
And in terms of the mapping $F$ defined by \eqref{mapping}, the strong version of the Eshelby conjecture can be mathematically interpreted as
\beqs
\begin{split}
&\forall\;(\boldsymbol{\sigma}^*, \bfC)\in \{\boldsymbol{\sigma ^*}\}\times \{\bfC\},\;\mathrm{if} \;\;\exists\; \Omega\in \{\Omega\} \;\;\mathrm{such\;that} \;F(\boldsymbol{\sigma^*},\;\bfC,\;\Omega)=\mathrm{constant},\;\mathrm{then}\;\Omega\in\{E\}.
\end{split}
\eeqs

 We let $\overline{\boldsymbol{\sigma}}^{*}$ denote a special kind of the uniform eigenstress that possesses two identical eigenvalues $k_1$ and a distinct eigenvalue $k_3$ with the same sign as that of $k_1$. And we let  $\{\overline{\boldsymbol{\sigma}}^{*}\}$ denote the set of $\overline{\boldsymbol{\sigma}}^{*}$, which is a subset of $\{{\boldsymbol{\sigma}}^{*}\}$.
 Then Theorem \ref{T0.3} can be mathematically interpreted as
 \beqs\label{T0.3math}
 \begin{split}
 \;&\forall\;\overline{\boldsymbol{\sigma}}^{*}\in \{\overline{\boldsymbol{\sigma}}^{*}\}, \;\;\exists\; \overline{\bfC}\in\{\bfC\},\;\Omega\in \{\Omega\},\; E\in\{E\} \;\mathrm{with}\;E\supset\Omega,\\
\; & \mathrm{such\;that} \;\; F(\overline{\boldsymbol{\sigma}}^{*},\;\overline{\bfC},\;\Omega)=F(\overline{\boldsymbol{\sigma}}^{*},\;\overline{\bfC},\;E)\equiv \mathrm{constant},\;\text{if and only if}\;\Omega\in\{E\}.
\end{split}
 \eeqs
 %\beqs\label{T0.3math}
% \begin{split}
% \;&\forall\;\overline{\boldsymbol{\sigma}}^{*}\in \{\overline{\boldsymbol{\sigma}}^{*}\},\;\mathrm{if} \;\;\exists\; \Omega\in \{\Omega\},\; E\in\{E\} \;\mathrm{with}\;E\supset\Omega,\\
 %\;&\mathrm{and \; meanwhile},\;\; \exists\; f:\{\overline{\boldsymbol{\sigma}}^{*}\}\rightarrow \{\overline{\bfC}\}\;\mathrm{with}\;\{\overline{\bfC}\}\subset \{\bfC\},\\
%\; & \mathrm{such\;that} \;\; %F(\overline{\boldsymbol{\sigma}}^{*},f(\overline{\boldsymbol{\sigma}}^{*}),\Omega)=F(\overline{\boldsymbol{\sigma}}^{*},f(\overline{\boldsymbol{\sigma}}^{*}), E)\equiv \mathrm{constant},\\
%\; &\mathrm{then}\;\Omega\in\{E\},
% \end{split}
% \eeqs

 We use $\{\overline{\bfC}\}$ to denote the set of $\overline{\bfC}$ that makes \eqref{T0.3math} valid.
Then we can conclude from \eqref{T0.3math} that for the uniform elastic strain fields  constrained to that induced by an ellipsoid, i.e.,  $F(\overline{\boldsymbol{\sigma}}^{*},\;\overline{\bfC},\;E)$,
from the specified uniform eigenstress $\overline{\boldsymbol{\sigma}}^{*}$, the strong version of the Eshelby conjecture is true for a set $\{\overline{\bfC}\}$ of isotropic elastic tensors, which is associated with the set $\{\overline{\boldsymbol{\sigma}}^{*}\}$ of the specified uniform eigenstresses.
In this regard, we will show that the Lam\'{e} parameters $\lambda,\mu$ for the isotropic elastic tensor $\overline{\bfC}$ are determined by the eigenvalues $k_1,k_3$ of the specified uniform eigenstress $\overline{\boldsymbol{\sigma}}^{*}$.

Then we turn to prove these three theorems.

\subsection{ Proof of Theorem \ref{T0.2}}

It can be directly derived from \eqref{NandNtilde} that
\begin{align}\label{algebraic}
  \begin{split}
  N_{\Omega }(\bfx)&=\frac{q_2(\bfx)+\varphi_2(x_1,x_2)-q_1(\bfx)-\varphi_1(x_3)}{\eta-\gamma},\\
 \tilde{N}_{\Omega }(\bfx) &= \frac{\eta q_1(\bfx)+\eta\varphi_1(x_3)-\gamma q_2(\bfx)-\gamma\varphi_2(x_1,x_2)}{\eta-\gamma}\quad\bfx\in \Omega.
  \end{split}
  \end{align}
By substituting $\eqref{algebraic}_1$ into \eqref{definition N}, we gain
\begin{align}
 \frac{\partial^2 \varphi_1(x_3)}{\partial x_3^2}+\frac{\partial^2 \varphi_2(x_1,x_2)}{\partial x_1^2}+\frac{\partial^2 \varphi_2(x_1,x_2)}{\partial x_2^2} =\mathrm{constant},
\end{align}
which indicates that $\varphi_1(x_3)$ is a constant, linear or quadratic function, and $\varphi_2(x_1,x_2)$ must satisfy
\begin{align}
\Delta \varphi_2(x_1,x_2)=\mathrm{constant}.
\end{align}

Up to now, it is straightforward to see that $\eqref{algebraic}_1 $ can be rewritten as
\begin{align}\label{Newtonianplannar}
{N_{\Omega }}(\bfx)=q'(\bfx)+\varphi'(x_1,x_2)\quad\bfx'\in\Omega,
\end{align}
where $q'(\bfx):=\frac{q_2(\bfx)-q_1(\bfx)-\varphi_1(x_3)}{\eta-\gamma}$, and $\varphi'(x_1,x_2):=\frac{\varphi_2(x_1,x_2)}{\eta-\gamma}$.

Then we are going to derive \eqref{Ncorrelation} from \eqref{Newtonianplannar}.
 Since $\Omega$ is assumed to be convex, $\forall \;\bfp \in\partial \Omega$, there must be a supporting hyperplane $S^*(\bfp)$ that passes through $\bfp$ but does not intersect the interior of $\Omega$, which means $\Omega$ is entirely on one side of $S^*(\bfp)$. Hence, if we define $\bfn_{p}$ as the normal vector of $S^*(\bfp)$, which points outside $\Omega$, we see that
\begin{align}\label{p}
\begin{split}
\bfn_{p}\cdot(\bfp-\bfy)>0,
\;\;\forall \bfy\in \Omega\;\;\text{and}\;\;\bfy \neq \bfp,
\end{split}
\end{align}
and thus
\begin{align}\label{minimun}
\begin{split}
\left. {{\bfn_{{p}}} \cdot \;{\boldsymbol{\nabla} }_\bfx{N_{\Omega}}(\bfx)} \right|_{\bfp} =  \int _{\Omega} \frac{\bfn_{p}\cdot(\bfp-\bfy)}{{4\pi {{\left| \bfp-\bfy \right|}^3}}}d\bfy  > 0,
\end{split}
\end{align}
which indicates that ${N_{\Omega }}(\bfx)$ always increases when $\bfx$ moves from the interior of $\Omega$ to the exterior of $\Omega$. Hence there must be a minimum point $\bfM$ of ${N_{\Omega }}\left( \bfx \right)$ inside $\Omega$. Since ${N_{\Omega }}(\bfx)$ is analytic inside $\Omega$,  we can expand ${N_{\Omega }}(\bfx)$ into the Taylor series at $\bfM$ and then take $\bfM$ as the origin of a new Cartesian coordinate system $\bfx'=(x_1',x_2',x_3')$, which yields
\begin{align}\label{ha}
\begin{split}
{N_{\Omega }}(\bfx')&=-\int_\Omega\frac{1}{4\pi|\bfx'-\bfy'|}d\bfy'\\
&={N_{\Omega }}(0)+\bfx'\cdot\nabla\otimes\nabla{N_{\Omega }}(0)\cdot\bfx'+\varphi^0(\bfx')\;\;\bfx'\in U(0),
\end{split}
\end{align}
where $\bfy'=(y_1',y_2',y_3')$ are the coordinates of the point within $\Omega$ on the basis of the new coordinate system; $\varphi^0(\bfx')$ denotes the sum of the terms whose degree is larger than 2; and $U(0)\subset \Omega$ is the neighborhood of the origin. Note that
$\nabla\otimes\nabla{N_{\Omega }}(0)$ is positive definite due to the origin being the minimum point of $N_{\Omega }$.

It can be seen from \eqref{Newtonianplannar} that in the new coordinate system $\bfx'=(x_1',x_2',x_3')$,
\begin{align}\label{haha}
{N_{\Omega }}(\bfx')=q'(\bfx')+\varphi'(x_1',x_2')\quad\bfx'\in\Omega.
\end{align}
   Then we can expand the right-hand side of \eqref{haha} into the Taylor series at the origin, which results in
\begin{align}\label{hahaha}
{N_{\Omega }}(\bfx')=q^*(\bfx')+\varphi^*(x_1',x_2')\quad\bfx'\in U(0),
\end{align}
where $q^*(\bfx')$ is a quadratic function, and $\varphi^*(x_1',x_2')$ is  the sum of the terms with respect to $x_1'$ and $x_2'$ whose degree is larger than 2. %It is noted {\color{red} (how?!)} that there must be a quadratic function

Then comparison of \eqref{hahaha} with \eqref{ha}  leads to
\begin{align}\label{haha7}
\varphi^*(x_1',x_2')=\varphi^0(\bfx')
\end{align}
 and
 \begin{align}\label{haha777}
 q^*(\bfx')={N_{\Omega }}(0)+\bfx'\cdot\nabla\otimes\nabla{N_{\Omega }}(0)\cdot\bfx'.
   \end{align}
%And since \eqref{haha} is also valid when $\bfx'\in U(0)$, substituting \eqref{haha} into \eqref{hahaha} yields
%subtracting both sides of \eqref{haha} from \eqref{hahaha} simultaneously yields
Let
\begin{align}\label{hahaha777}
q^0(x_1',x_2'):=q'(\bfx')-q^*(\bfx');
\end{align}
%where $q^0(x_1',x_2')$ utilized here exactly represents a portion of the expansion of $\varphi'(x_1',x_2')$ at the origin.
thus by substituting \eqref{haha777} and \eqref{hahaha777} into \eqref{haha}, we obtain
\begin{align}\label{ultimate1}
\begin{split}
{N_{\Omega }}(\bfx')=\;&{N_{\Omega }}(0)+\bfx'\cdot\nabla\otimes\nabla{N_{\Omega }}(0)\cdot\bfx'+q^0(x_1',x_2')+\varphi'(x_1',x_2')\quad\quad\quad\quad\bfx'\in\Omega.
\end{split}
\end{align}
Let $\varphi''(x_1',x_2'):=q^0(x_1',x_2')+\varphi'(x_1',x_2')$; it follows from \eqref{ultimate1} that
\begin{align}\label{ultimate2}
\begin{split}
&{N_{\Omega }}(\bfx')\!=\!{N_{\Omega }}(0)\!+\!\bfx'\cdot\nabla\otimes\nabla{N_{\Omega }}(0)\cdot\bfx'\!+\!\varphi''(x_1',x_2')\;\;\;\bfx'\in\Omega.
\end{split}
\end{align}
Note that the term $\bfx'\cdot\nabla\otimes\nabla{N_{\Omega }}(0)\cdot\bfx'$ on the right-hand side of \eqref{ultimate2} needs to satisfy
\begin{align}\label{constrainforN}
\Delta(\bfx'\cdot\nabla\otimes\nabla{N_{\Omega }}(0)\cdot\bfx')=1,
\end{align}
which is derived by substitution of \eqref{ha} into $\eqref{definition N}_1$.
Besides, we have already proved that $\bfx'\cdot\nabla\otimes\nabla{N_{\Omega }}(0)\cdot\bfx'$ is positive definite.

Then based on the above property of $\bfx'\cdot\nabla\otimes\nabla{N_{\Omega }}(0)\cdot\bfx'$, it is known from \cite{Benedetto1986} that there must be an ellipsoid $E\subset\Omega$ satisfying
\begin{align}\label{Neexpression}
\begin{split}
N_{E}(\bfx')=C^E+\bfx'\cdot\nabla\otimes\nabla{N_{\Omega }}(0)\cdot\bfx'\quad\bfx'\in E
\end{split}
\end{align}
with $C^E$ a positive real constant.

Then substituting \eqref{Neexpression} into \eqref{ultimate2} yields
\begin{align}
\begin{split}
N_{\Omega}(\bfx')= N_{E}(\bfx')+{\varphi}(x_1',x_2')\quad\bfx'\in E,
\end{split}
\end{align}
where $\varphi(x_1',x_2'):=\varphi''(x_1',x_2')+{N_{\Omega }}(0)-C^E$. Thus the proof of Theorem \ref{T0.2} is completed.

\subsection{ Proof of Theorem \ref{new}}

 First of all, because we are going to consider two elastic tensors, denoted by $\bfC^{(1)}$ and $\bfC^{(2)}$,  all of the parameters, variables, and functions that correspond to $\bfC^{(1)}$ and $\bfC^{(2)}$ will be distinguished by the superscripts (1) and (2), respectively.

 Since \eqref{NandNtilde} holds for $( \boldsymbol{\sigma^*},\bfC^{(1)})$ and $( \boldsymbol{\sigma^*},\bfC^{(2)})$ simultaneously,
  it is derived from $\eqref{NandNtilde}_1$ that for $\bfx\in \Omega$,
 \begin{align}\label{algebraicnew}
  \begin{split}
  \gamma^{(i)} N_\Omega(\bfx)+\tilde{N}_{\Omega}(\bfx)=q^{(i)}_1(\bfx)+\varphi^{(i)}_1(x_3)\quad i=1,2,
  \end{split}
  \end{align}
  where
  \begin{align}
  \gamma^{(i)}=\frac{[k_3(\lambda^{(i)}+\mu^{(i)})-k_1(\lambda^{(i)}+3\mu^{(i)})]}{(\lambda^{(i)}+\mu^{(i)})(k_1-k_3)}\quad i=1,2;
    \end{align}
  $q^{(i)}_1(\bfx)\;(i=1,2)$ are quadratic functions; and $\varphi^{(i)}_1(x_3)\;(i=1,2)$ have been proved to be constant, linear or quadratic functions in the above derivation.

  According to Theorem \ref{T0.2}, if $\Omega$ is convex, there exists an ellipsoid $E$ such that \eqref{Ncorrelation} holds. Thus by substituting \eqref{Ncorrelation} into \eqref{algebraicnew}, we gain
  \begin{align}\label{phi2}
    \begin{split}
\varphi_2(x_1,x_2)=&\frac{q^{(1)}_1(\bfx)+\varphi^{(1)}_1(x_3)-[q^{(2)}_1(\bfx)+\varphi^{(2)}_1(x_3)]}{\gamma^{(1)}-\gamma^{(2)}}\\
&-N_E(\bfx)\quad\quad\quad\quad\quad\quad\quad\quad\quad\quad\bfx\in\Omega.
  \end{split}
   \end{align}
   The validity of \eqref{phi2} calls for
    \begin{align}
    \gamma^{(1)}\neq\gamma^{(2)},
  \end{align}
which requires
\begin{align}\label{constraint1}
\lambda^{(1)}\mu^{(2)}-\lambda^{(2)}\mu^{(1)}\neq 0.
\end{align}
Note that \eqref{constraint1} is satisfied owing to the linear independence of $\bfC^{(1)}$ and $\bfC^{(2)}$.

 Since the Newtonian potential $N_E(\bfx)$ of the ellipsoid $E$ is a quadratic function of $\bfx$ \cite{Ferrers,Dyson}, it is concluded from \eqref{phi2} that $\varphi_2(x_1,x_2)$ can only be a constant, linear or quadratic function of $x_1,x_2$.

Given that $\varphi_2(x_1,x_2)$ can only be a constant, linear or quadratic function, we see that $N_\Omega(\bfx)$ must be quadratic inside $\Omega$ due to \eqref{Ncorrelation}. To further determine the shape of $\Omega$, we resort to a powerful theorem, i.e.,
\begin{theorem}\label{TheoremNP}
 Let $\Omega$ be a bounded domain with a Lipschitz boundary. The relation
\begin{align}
 N_\Omega(\bfx)=-\frac{1}{4\pi}\int_{\Omega} \frac{1}{{|\bfx-\bfy|}} \mathrm{d}\bfy =\mathrm{quadratic}\quad \bfx\in\Omega
 \end{align}
holds if and only if $\Omega$ is an ellipsoid~\cite{Kang2008}.
\end{theorem}
 %Dive \cite{Dive1931} and Nikliborc \cite{Nikliborc1932} first substantiated Theorem~\ref{TheoremNP} for $C^1$ domains, and then Kang and %Milton\cite{Kang2008} ultimately substantiated Theorem~\ref{TheoremNP} for Lipschitz domains.
Since we have proved that $N_\Omega(\bfx)$ must be quadratic inside $\Omega$, then according to Theorem~\ref{TheoremNP}, we conclude that $\Omega$ must be ellipsoidal. Thus  Theorem \ref{new} is proved for convex inclusions.

 We note that even if the inclusion $\Omega$ is concave, we can also verify Theorem \ref{new}. It can be directly derived  from \eqref{algebraicnew} that
  \begin{align}\label{Newtoniannew}
    \begin{split}
    N_\Omega(\bfx)=\frac{q^{(1)}_1(\bfx)+\varphi^{(1)}_1(x_3)-[q^{(2)}_1(\bfx)+\varphi^{(2)}_1(x_3)]}{\gamma^{(1)}-\gamma^{(2)}}\quad\bfx\in\Omega.
     \end{split}
  \end{align}
Since \eqref{Newtoniannew} still indicates that $N_\Omega(\bfx)$ must be quadratic inside $\Omega$, then by Theorem \ref{TheoremNP}, we draw the conclusion that $\Omega$ is of ellipsoid shape, which completes the proof of Theorem \ref{new}.

\subsection{Proof of Theorem \ref{T0.3}}

According to the derivations from \eqref{disfourier} to \eqref{NandNtilde},  Theorem \ref{T0.3} indicates that there exists an ellipsoid $E$ satisfying $E\supset\Omega$, such that
\begin{align}\label{newresult}
\begin{split}
 &\frac{\partial}{\partial x_{i}}\left[\gamma N_\Omega(\bfx)+\tilde{N}_{\Omega}(\bfx)\right]= \frac{\partial}{\partial x_{i}}\left[\gamma N_E(\bfx)+\tilde{N}_{E}(\bfx)\right]\quad i=1,2,\\
&\frac{\partial}{\partial x_3}\left[\eta N_\Omega(\bfx)+\tilde{N}_{\Omega}(\bfx)\right]=\frac{\partial}{\partial x_3}\left[\eta N_E(\bfx)+\tilde{N}_{E}(\bfx)\right]\quad\bfx\in \Omega,
 \end{split}
\end{align}
where $\tilde{N}_{E}(\bfx)$ is given by \eqref{Ntildeexpression} with $\Omega$ replaced by $E$.

Let $E^*=\{t\bfx\;|\;\bfx\in E, t>0\}$, which is still large enough to contain $\Omega$ but intersects with $\Omega$ at $\bfQ=(Q_1,Q_2,Q_3)$. Then owing to the quadratic forms of $N_{E}(\bfx)$ and $\tilde{N}_{E}(\bfx)$, which can be explicitly calculated for any $E$, it is easy to verify that
\begin{align}\label{quadraticNandNT}
N_{E^*}(\bfx)-N_{E}(\bfx)=\mathrm{constant},\;\;\tilde{N}_{E^*}(\bfx)-\tilde{N}_{E}(\bfx)=\mathrm{constant},
\end{align}
substitution of which into \eqref{newresult} yields
\begin{align}\label{original3}
  \begin{split}
&\frac{\partial}{\partial x_{i}}\left[\gamma N_{E^*\setminus\Omega }(\bfx)+\tilde{N}_{E^*\setminus\Omega }(\bfx)\right]= 0\quad i=1,2,\\
&\frac{\partial}{\partial x_3}\left[\eta N_{E^*\setminus\Omega }(\bfx)+\tilde{N}_{E^*\setminus\Omega }(\bfx)\right]=0\quad\bfx\in\Omega,
  \end{split}
\end{align}
where
 \begin{align}\label{Ntildeage}
  \begin{split}
  N_{E^*\setminus\Omega }(\bfx)&:=-\frac{1}{4\pi}\int_{E^*\setminus\Omega }\frac{1}{|\bfx-\bfy|}\mathrm{d}\bfy,\\
\tilde{N}_{E^*\setminus\Omega }(\bfx)&:=-\frac{1}{4\pi}\int_{E^*\setminus\Omega }\frac{(x_3-y_3)^2}{|\bfx-\bfy|^3}\mathrm{d}\bfy.
  \end{split}
  \end{align}

 Then we can choose an isotropic elastic tensor $\bfC$ that leads to
   \begin{align}\label{choose0}
   \gamma=0,
   \end{align}
   which requires
  \begin{align}\label{choose1}
k_1(\lambda+3\mu)=k_3(\lambda+\mu).
\end{align}
 Under the condition that $k_1$ and $k_3$ have the same sign, it is easy to verify that there exist $\lambda,\mu$ that satisfy \eqref{Lam} such that \eqref{choose1} holds. Note that if $k_1$ and $k_3$ have opposite signs, \eqref{choose1} cannot be valid.

 Then in terms of \eqref{choose0},  it can be derived from \eqref{original3} that
  \begin{align}\label{algebraic2}
  \begin{split}
  N_{E^*\setminus\Omega }(\bfx)&=\frac{1}{\eta}[\varphi_2(x_1,x_2)+C_2-\varphi_1(x_3)-C_1],\\
  \tilde{N}_{E^*\setminus\Omega }(\bfx)&=C_1+\varphi_1(x_3)\quad\quad\bfx\in \Omega,
  \end{split}
  \end{align}
  where $C_1,C_2$ denote two real constants, and $\varphi_1,\varphi_2$ still denote two unknown functions.

  Let $\bfn=(n_1,n_2,n_3)$ denote the outward-pointing normal vector of $\partial E^*$ at the point $\bfQ=(Q_1,Q_2,Q_3)$.  Then it can be derived from $\eqref{Ntildeage}_2$ that
  \begin{align}\label{derivative}
  \begin{split}
  \left.\bfn\cdot \boldsymbol{\nabla}\tilde{N}_{E^*\setminus\Omega }\right|_\bfQ
  =&-2n_3\int_{E^*\setminus\Omega }\frac{(Q_3-y_3)}{4\pi|\bfQ-\bfy|^3}\mathrm{d}\bfy+3\int_{E^*\setminus\Omega }\frac{(Q_3-y_3)^2(\bfQ-\bfy)\cdot\bfn}{4\pi|\bfQ-\bfy|^5}\mathrm{d}\bfy\\
  =&\left.-2n_3\frac{\partial N_{E^*\setminus\Omega }}{\partial x_3}\right|_\bfQ+3\int_{E^*\setminus\Omega }\frac{(Q_3-y_3)^2(\bfQ-\bfy)\cdot\bfn}{4\pi|\bfQ-\bfy|^5}\mathrm{d}\bfy.
  \end{split}
  \end{align}

We define  a vector function $\bfF$, i.e.,
\begin{align}\label{definitionF}
\bfF(\bfx):=3\int_{E^*\setminus\Omega }\frac{(x_3-y_3)^2(\bfx-\bfy)}{4\pi|\bfx-\bfy|^5}\mathrm{d}\bfy,
  \end{align}
  and thus substituting \eqref{definitionF} into \eqref{derivative} leads to
  \begin{align}\label{makeargument}
 \left.\bfn\cdot \bfF\right|_\bfQ= \left.\left(\bfn\cdot \boldsymbol{\nabla}\tilde{N}_{E^*\setminus\Omega }+2n_3\frac{\partial N_{E^*\setminus\Omega }}{\partial x_3}\right)\right|_\bfQ.
   \end{align}

   Based on \eqref{makeargument}, we will make some arguments on the change of $\bfF$ at $\bfQ$ near the boundary $\partial E^*$, just like the arguments on a similar situation for the change of the Newtonian potential made by Kang and Milton \cite{Kang2008}.

Likewise, we consider two cases concerning the continuity of $\partial\Omega$.

\subsubsection{ $\partial\Omega$ possesses $C^1$ continuity}

  In this case, $\bfn$ is also the outward-pointing normal vector of  $\partial \Omega$ at $\bfQ$.
  Thus by substituting \eqref{algebraic2} into \eqref{makeargument}, we obtain
  \begin{align}\label{3}
 \left.\bfn\cdot \bfF\right|_\bfQ= \left.(\frac{2}{{\eta }}-1)n_3\frac{\partial \varphi_1(x_3)}{\partial x_3}\right|_\bfQ.
\end{align}
{To analyze the change of $\bfF$ at $\bfQ$ via \eqref{3}, we eliminate the unknown function $\varphi_1(x_3)$ on the right-hand side of \eqref{3} by letting
 \begin{align}\label{c1condition}
\eta=2,
\end{align}
which requires
 \begin{align}\label{choose2}
k_1(\lambda+\mu)=k_3(\mu-\lambda).
\end{align}
%Recall that for given eigenvalues $k_1$ and $k_3$, we have chosen $\lambda,\mu$ that satisfy \eqref{choose1}.
Note that for given eigenvalues $k_1$ and $k_3$, \eqref{choose1} and \eqref{choose2} together determine a unique combination $(\lambda,\mu)$ of the Lam\'{e} parameters  and thus a unique isotropic elastic tensor $\bfC$.
%From another perspective, \eqref{choose1} and \eqref{choose2} constitute the relationship between defined above.

%Since $\lambda,\mu$ are two mutually independent parameters, {\color{blue} under the condition \eqref{choose1},} we can further choose $\lambda,\mu$ that yield
 %We note that for specific $k_1$ and $k_3$, \eqref{choose1} and \eqref{choose2} determine a unique isotropic elastic tensor $\bfC$. Or in other words, the corresponding Lam\'{e} parameters $\lambda,\mu$ that satisfy \eqref{Lam} are determined by \eqref{choose1} and \eqref{choose2} when $k_1$ and $k_3$ are given. {\color{red}\bf (it seems that (\ref{choose2}) is artificial and quite arbitrary.....which is simply only an example? i think
%(\ref{choose1}) is the only condition that determines the set of
%C! right? if so, we must take (\ref{choose2}) as an example!)} {\color{blue} (For given $k_1$ and $k_3$, we do not choose a kind of isotropic media (a set of $\bfC$) to verify the third theorem, but choose a unique $\bfC$ to verify the third theorem. \eqref{choose1} and \eqref{choose2} together determines a unique $\bfC$ ($\lambda$ and $\mu$) for given $k_1$ and $k_3$, and only in this unique isotropic medium, our argument can finally work.) }
Then by substituting \eqref{c1condition} into \eqref{3}, we obtain
  \begin{align}\label{newcontradiction}
 \left.\bfn\cdot \bfF\right|_\bfQ= 0.
\end{align}
However, if $E^*\setminus\Omega $ is not empty,  we see $(\bfQ-\bfy)\cdot\bfn\geq 0$, substitution of which into \eqref{definitionF} yields
\begin{align}\label{2}
 \left.\bfn\cdot \bfF\right|_\bfQ=3\int_{E^*\setminus\Omega }\frac{(Q_3-y_3)^2(\bfQ-\bfy)\cdot\bfn}{4\pi|\bfQ-\bfy|^5}\mathrm{d}\bfy>0.
\end{align}
Since \eqref{2} contradicts \eqref{newcontradiction}, thus $E^*\setminus\Omega$ must be empty to avoid the contradiction, which leads to the conclusion $\Omega=E^*$.

  Note that \eqref{3} is valid only when there are line segments $\{\bfQ\pm t\bfn\;|t\in \mathbb{R}\}$ belonging to $\Omega$, which is obvious if $\partial\Omega$ is $C^1$ continuous.

\subsubsection{$\partial\Omega$ possesses Lipschitz continuity}

  In this case, as is stated before, \eqref{3} may not hold since there may not be a line segment $\{\bfQ\pm t\bfn\;|t\in \mathbb{R}\}$  belonging to $\Omega$.
  Hence there is an alternative way to proceed with the argument.

  If $E^*\setminus\Omega $ is not empty, from \eqref{2}, we know that the vector function $\bfF$ at $\bfQ$ must satisfy
  \begin{align}\label{7}
  \left.\bfF\right|_\bfQ\neq \bf{0}.
  \end{align}
  Let $L$ be a line passing through $\bfQ$ and satisfying $L\cap\Omega\subset\Omega$. The direction of $L$ is represented by $\bfv$. Since $\partial\Omega$ has Lipschitz continuity, there is a neighborhood $\bfV$ in $S^2$ of $\bfv$, where $S^2$ denotes the unit sphere  in $\mathbb{R}^3$; any line $L^0$ passing through $\bfQ$ and possessing the direction of the vector in $\bfV$ will lead to $L^0\cap\Omega\subset\Omega$. Those lines are contained in a set denoted by $\{L^0\}$. Based on \eqref{7}, we know that $\exists \;L^*\in \{L^0\}$, whose direction is $\bfv^*$, such that
  \begin{align}\label{4}
  \bfv^*\cdot\left.\bfF\right|_\bfQ\neq 0.
  \end{align}
  However, since $L^*\in \{L^0\}$ so that $L^*\cap\Omega\subset\Omega$,  there are line segments $\{\bfQ\pm t\bfv^*\;|t\in \mathbb{R}\}$ belonging to $\Omega$, and thus \eqref{3} can still be valid with $\bfn$ replaced by $\bfv^*$, which leads to
  \begin{align}\label{c2condition}
 \left.\bfv^*\cdot \bfF\right|_\bfQ= \left.(\frac{2}{{\eta }}-1)v^*_3\frac{\partial \varphi_1(x_3)}{\partial x_3}\right|_\bfQ.
\end{align}
Finally, by similarly letting \eqref{c1condition} and then
substituting \eqref{c1condition} into \eqref{c2condition}, we obtain
   \begin{align}
  \bfv^*\cdot\left.\bfF\right|_\bfQ= 0,
  \end{align}
  which contradicts \eqref{4}. Likewise, $E^*\setminus\Omega$ must be empty to avoid the contradiction so that $\Omega=E^*$. Therefore, the proof of Theorem \ref{T0.3} is fulfilled.

\section{Conclusions}

In this work, we have studied the strong version of the Eshelby conjecture in the context of three-dimensional isotropic elasticity, for the case where two of the eigenvalues of the eigenstress are identical and the other one is distinct.  We have made progress towards the proof of the conjecture  by presenting and proving three theorems. The first
theorem gives a necessary condition for the Newtonian potential of convex inclusions possessing the Eshelby uniformity property, which can exclude  non-ellipsoidal convex inclusions that do not satisfy the condition.  In terms of the necessary condition, and by introducing the effect of the elastic constants of the isotropic medium, the second theorem indicates that
there exist combinations of the elastic tensors and uniform eigenstresses
such that only an ellipsoid can have the Eshelby uniformity property for these combinations simultaneously.
The third theorem provides a more specifically constrained proof of the conjecture. It proves that
for the uniform strain fields constrained to that induced by an ellipsoid
from a set of the specified uniform eigenstress, the strong version of the Eshelby conjecture is true for a set of isotropic elastic tensors which are associated with the specified uniform eigenstress.

\begin{flushleft}
\textbf{\emph{Acknowledgements}}
\end{flushleft}
The authors sincerely thank Professor Liping Liu of Rutgers University for helpful discussions. The authors also acknowledge the support of the National Natural Science Foundation of China under Grant 11521202.

\appendix
\setcounter{equation}{0}
\renewcommand\theequation{A.\arabic{equation}}

\section*{Appendix}
\section*{Appendix. A proof of the strong version for the eigenstress possessing either all distinct or all identical eigenvalues by using Fourier analysis}

We are going to prove the strong version for the eigenstress possessing either all distinct or all identical eigenvalues through an alternative method, i.e., the Fourier analysis, which is somehow more concise than the method proposed by Ammari et al. \cite{Ammari2010}.

We then consider the two cases concerning the eigenvalues, separately.

\subsubsection*{(1) $\boldsymbol{\sigma}^*$ has three identical eigenvalues}

Assume $\Omega$ possesses the Eshelby uniformity property; thus \eqref{finaluineigenspace} holds.
Let $\sigma^*_{11}=\sigma^*_{22}=\sigma^*_{33}=k$. Then by substituting $\sigma^*_{11}=\sigma^*_{22}=\sigma^*_{33}=k$ along with \eqref{Green iso} and  \eqref{gsimplification}  into \eqref{finaluineigenspace} and then substituting \eqref{finaluineigenspace} into \eqref{NewtonianFourier}, we obtain
\begin{align}\label{identicalcondition}
\nabla\otimes\bfu=\frac{k}{\lambda+2\mu}\nabla\otimes\nabla N_\Omega(\bfx)=\mathrm{constant}\quad\bfx\in \Omega,
\end{align}
which indicates that $N_\Omega(\bfx)$ must be quadratic inside $\Omega$.  Then by resorting to Theorem~\ref{TheoremNP}, we conclude that $\Omega$ must be ellipsoidal.

\subsubsection*{(2) $\boldsymbol{\sigma}^*$ has three distinct eigenvalues}

Likewise, assume $\Omega$ possesses the Eshelby uniformity property; thus \eqref{finaluineigenspace} holds. Let $\sigma^*_{11}=k_1, \sigma^*_{22}=k_2,\sigma^*_{33}=k_3$ with $k_1\neq k_2,\;k_2\neq k_3$,\;$k_3\neq k_1$.
Then by substituting $\sigma^*_{11}=k_1, \sigma^*_{22}=k_2,\sigma^*_{33}=k_3$ along with \eqref{Green iso} and  \eqref{gsimplification} into \eqref{finaluineigenspace}, we obtain

\begin{align}\label{distinct}
 \begin{split}
 \forall j=1,2,3,\quad&\int_{\mathbb{R}^3}\!\!\! \frac{\xi_1\xi_j}{\mu} \left(\frac{ k_1}{|\boldsymbol{\xi}|^2}-\frac{\lambda+\mu}{\lambda+2\mu}\frac{k_1\xi_1^2+k_2\xi_2^2+k_3\xi_3^2}{|\boldsymbol{\xi}|^4}\right)g(\bfx,\boldsymbol{\xi})\mathrm{d}\boldsymbol{\xi}=\mathrm{constant};\\
& \int_{\mathbb{R}^3} \!\!\!\frac{\xi_2\xi_j}{\mu} \left(\frac{ k_2}{|\boldsymbol{\xi}|^2}-\frac{\lambda+\mu}{\lambda+2\mu}\frac{k_1\xi_1^2+k_2\xi_2^2+k_3\xi_3^2}{|\boldsymbol{\xi}|^4}\right)g(\bfx,\boldsymbol{\xi})\mathrm{d}\boldsymbol{\xi}=\mathrm{constant};\\
 &\int_{\mathbb{R}^3}\!\!\! \frac{\xi_3\xi_j}{\mu} \left(\frac{ k_3}{|\boldsymbol{\xi}|^2}-\frac{\lambda+\mu}{\lambda+2\mu}\frac{k_1\xi_1^2+k_2\xi_2^2+k_3\xi_3^2}{|\boldsymbol{\xi}|^4}\right)g(\bfx,\boldsymbol{\xi})\mathrm{d}\boldsymbol{\xi}=\mathrm{constant}\quad\bfx\in \Omega.\\
 \end{split}
\end{align}
 
It can be derived from \eqref{distinct} that
\begin{align}\label{distinctcondition}
 \begin{split}
&\int_{\mathbb{R}^3}\frac{\xi_1\xi_2 (k_1-k_2)}{|\boldsymbol{\xi}|^2}g(\bfx,\boldsymbol{\xi})\mathrm{d}\boldsymbol{\xi}=\mathrm{constant};\\
&\int_{\mathbb{R}^3}\frac{\xi_2\xi_3 (k_2-k_3)}{|\boldsymbol{\xi}|^2}g(\bfx,\boldsymbol{\xi})\mathrm{d}\boldsymbol{\xi}=\mathrm{constant};\\
&\int_{\mathbb{R}^3}\frac{\xi_3\xi_1 (k_3-k_1)}{|\boldsymbol{\xi}|^2}g(\bfx,\boldsymbol{\xi})\mathrm{d}\boldsymbol{\xi}=\mathrm{constant}\quad\bfx\in\Omega.
 \end{split}
\end{align}

We note that the integrals on the left-hand side of \eqref{distinctcondition} are exactly the second derivatives of the Newtonian potential $N_\Omega(\bfx)$, i.e.,
\begin{align}\label{Newtonianproperty}
\frac{\partial^2 N_\Omega(\bfx)}{\partial x_p\partial x_l}=\frac{1}{(2\pi)^3}\int_{\mathbb{R}^3} \frac{\xi_j\xi_l}{|\boldsymbol{\xi}|^2}g(\bfx,\boldsymbol{\xi})\mathrm{d}\boldsymbol{\xi},
\end{align}
which can be derived from \eqref{NewtonianFourier}.

Since $k_1\neq k_2,\;k_2\neq k_3$,\;$k_3\neq k_1$, by comparing \eqref{distinctcondition} with \eqref{Newtonianproperty}, we see
\begin{align}
 \forall i\neq j,\quad\frac{\partial^2 N_\Omega(\bfx)}{\partial x_i\partial x_j}=\mathrm{constant}\quad\bfx\in\Omega,
\end{align}
which implies
\begin{align}\label{Ndistinct}
 N_\Omega(\bfx)=q(\bfx)+\psi_1(x_1)+\psi_2(x_2)+\psi_3(x_3)\quad\bfx\in \Omega,
\end{align}
 where $q$ denotes a quadratic function, and $\psi_i(i=1,2,3)$ denote unknown functions.

 By substituting \eqref{Ndistinct} into $\eqref{definition N}_1$, we gain
\begin{align}
 \frac{\partial^2\psi_1(x_1)}{\partial x_1^2}+\frac{\partial^2\psi_2(x_2)}{\partial x_2^2}+\frac{\partial^2\psi_3(x_3)}{\partial x_3^2}=\mathrm{constant}.
\end{align}
Hence $\psi_i\ (i=1,2,3)$ must be constant, linear or quadratic functions, and thus the Newtonian potential $N_\Omega(\bfx)$ in \eqref{Ndistinct} must be quadratic. Likewise, according to Theorem~\ref{TheoremNP}, we conclude that $\Omega$ must be ellipsoidal.
%Hence we have proved the strong version of the Eshelby conjecture for the eigenstress possessing three distinct eigenvalues. Therefore, the proof of %the strong version of the Eshelby conjecture for the eigenstress possessing either all distinct or all identical eigenvalues is completed.
%Therefore, we draw the conclusion that in the three-dimensional isotropic medium, the ellipsoidal inclusion uniquely transforms a single uniform %eigenstress that has either all distinct or all identical eigenvalues into a uniform stress field in it.

%%%%%%%%%%%%%%%%%%%%%%%%%%%%%%%%%%%%%
%Reference
%\bibliographystyle{spmpsci}
%\bibliography{Eshelbyconjecturedatabase}

%\end{sloppypar}
\end{document}